\documentclass[12pt]{article}
\usepackage[centertags]{amsmath}
\usepackage{amsfonts}
\usepackage{amssymb}
\usepackage{amsthm}
\usepackage{amsmath,mathrsfs}
\usepackage{amsmath,amscd}
\usepackage{mathtools}
\usepackage[matrix,arrow,curve,cmtip]{xy}
\title{\textbf{Ground States of the $\infty$-categorical Grothendieck Construction}}
\author{Renaud Gauthier \footnote{2020 Math. Subj. Class: 18A25, 18D30, 18N40, 18N60, 68P05, 81S07. Keywords: Grothendieck construction, $\infty$-categories.} \\ \\}
\theoremstyle{definition}
\newtheorem{conj}{Conjecture}[subsection]

\newtheorem*{acknowledgments}{Acknowledgments}

\DeclareMathOperator*{\adj}{\rlarr}

\newcommand{\beq}{\begin{equation}}
\newcommand{\eeq}{\end{equation}}

\newcommand{\hrarr}{\hookrightarrow}
\newcommand{\rarr}{\rightarrow}

\newcommand{\rlarr}{\rightleftarrows}

\newcommand{\Ob}{\text{Ob\,}}

\newcommand{\xrarr}{\xrightarrow}

\newcommand{\cC}{\mathcal{C}}
\newcommand{\cCop}{\cC^{\op}}
\newcommand{\cD}{\mathcal{D}}
\newcommand{\cE}{\mathcal{E}}

\newcommand{\cU}{\mathcal{U}}

\newcommand{\cX}{\mathcal{X}}

\newcommand{\gC}{\mathfrak{C}}

\newcommand{\Cat}{\text{Cat}}

\newcommand{\Fun}{\text{Fun}}

\newcommand{\Hom}{\text{Hom}}

\newcommand{\op}{\text{op}}

\newcommand{\Set}{\text{Set}}

\newcommand{\Top}{\text{Top}}

\newcommand{\Catinf}{\Cat_{\infty}}

\newcommand{\CatD}{\Cat_{\Delta}}

\newcommand{\Dn}{\Delta^n}

\newcommand{\SetplusD}{\Set^+_{\Delta}}
\newcommand{\SetD}{\Set_{\Delta}}

\newcommand{\St}{\text{St}}

\newcommand{\Stphi}{\St_{\phi}}

\newcommand{\SetDplusK}{(\SetD^+)_{/K}}
\newcommand{\SetDplusC}{(\SetD^+)^{\cC}}

\newcommand{\Unphi}{\text{Un}_{\phi}}

\newcommand{\bp}{\textbf{p}}
\newcommand{\ED}{\cE_{\cD}}
\newcommand{\EX}{\cE_{X}}
\newcommand{\pBun}{\bp \text{-Bun}}
\newcommand{\pFib}{\bp \text{-Fib}}
\newcommand{\KaT}{K_{\alpha}^{\triangleleft}}
\newcommand{\pa}{p_{\alpha}}
\newcommand{\SetplusDp}{(\SetplusD)_{/\bp}}
\newcommand{\SetplusDSp}{(\SetplusD)_{/S \times \bp}}

\newcommand{\tgC}{\widetilde{\gC}}
\newcommand{\tgCS}{\widetilde{\gC}[S]}

\newcommand{\SetplusDptCS}{(\SetplusDp)^{\tgCS}}
\newcommand{\Unp}{\text{Un}\bp}
\newcommand{\SetplusDSUnp}{(\SetplusD)_{/S \times \Unp}}
\newcommand{\SetplusDS}{(\SetplusD)_{/S}}
\newcommand{\tSt}{\widetilde{\St}}
\newcommand{\tStplusphi}{\tSt^+_{\phi}}
\newcommand{\tStpplusphi}{\tSt^{'+}_{\phi}}
\newcommand{\tStphi}{\tSt_{\phi}}
\newcommand{\Azero}{A^{\circ}}
\newcommand{\gCK}{\gC[K]}
\newcommand{\coCFib}{\text{coCFib}}
\newcommand{\Stphiplus}{\Stphi^{+}}
\newcommand{\Unphiplus}{\Unphi^+}
\newcommand{\SetDplus}{\SetD^+}
\newcommand{\Xlhd}{X^{\lhd}}
\newcommand{\Un}{\text{Un}}
\newcommand{\tUnpphi}{\widetilde{\Un}^{\bp}_{\phi}}

\begin{document}
\maketitle
\begin{abstract}
We highlight the presence of a ground state in the $\infty$-categorical Grothendieck construction of Lurie, further developed by Arakawa, in which both straightened and unstraightened pictures coexist.
\end{abstract}

\newpage

\section{Introduction}
In \cite{RG} we saw that following a relative point of view philosophy, one can study a given geometric object $X$ by embedding it in a larger category $\cC$, an ambient category of some sort. This naturally gives rise to various categories, one of which is the slice category $\cC_{/X}$. Beyond that, one can either study the manifestation $\Hom(\cC_{/X},\cD)$ of $\cC_{/X}$ into some other category $\cD$, or one can study the perception of such a category $\cC_{/X}$ in the larger category it is an object of, a perception being constructed from maps into $\cC_{/X}$, as in $\Hom(\cE, \cC_{/X})$. The present paper has for aim to relate manifestations and perceptions to highlight their duality. This we do using a generalization of the Grothendieck construction in the $\infty$-categorical context. In particular we will use the straightening/unstraightening picture of \cite{Lu} as well as a particular generalization worked out in \cite{A}. In that regard, we will use those formalisms for the sake of developing ours, but regarding the $\infty$-categorical Grothendieck constructions themselves, we make no pretense to originality, the material can be found in \cite{Lu} and \cite{A}.\\

We will provide a phenomenological interpretation of the $\infty$-categorical Grothendieck construction, with a view towards applying this to Physical Mathematics. While doing so we will interpret the equivalence of $\infty$-categories $\coCFib(\cC) \simeq \Fun(\cC,\Catinf)$ as indicating that to probe the structure of $\cC$, which is provided by the left hand side of this equivalence, amounts to attaching $\infty$-categories to every object of $\cC$, as on the right hand side, thereby providing a microstructure to $\cC$. We will also introduce the notion of a ground state, a setting where dual models are present in a deconstructed form, at which point any duality has dissolved. Conversely, this means that the ground state can be regarded as a fundamental theory which condenses into two dual theories, a symmetry breaking of some sort. \\

Researchers that are new to this concept of Grothendieck construction in the $\infty$-categorical setting will find the present position helpful in shedding some light on some constructions, some of which are not exactly transparent, and researchers interested in the application of $\infty$-categorical techniques to Quantum Field Theory will find our exposition to be helpful. Within the same context we will also introduce the notion of a formal uncertainty principle.\\

Relation to past work: there is a related publication \cite{HB} within the philosophy of Science that defines a structure similar to our deconstruction of theories into a ground state, though in a more elementary setting, with a different objective, and with different implementations. Nevertheless, it is also worthy of consideration.\\

\section*{Notations}
We denote by $\SetD$ the category of simplicial sets, by $\CatD$ the category of simplicial categories, and by $\Catinf$ the $\infty$-category of $\infty$-categories. Given a model category $A$, $\Azero$ denotes its subcategory of fibrant and cofibrant objects. We refer to objects of the form $\Fun(\cC,\cD)$ as manifestations of $\cC$ into $\cD$, or as perceptions of $\cD$ from $\cC$. We will use the abbreviation RLP to mean ``Right Lifting Property".

\begin{acknowledgments}
The author would like to warmly thank Kensuke Arakawa for patiently explaining fundamental concepts of his paper.
\end{acknowledgments}

\section{Original Grothendieck construction}
We first give a lightning review of the Grothendieck construction (\cite{AG}). For $\cC$ a small category and $F: \cC \rarr \Cat$ a functor, one defines a category $\int F$ with objects pairs $(C,u)$, with $C \in \Ob(\cC)$ and $u \in \Ob(FC)$, and morphisms $(C,u) \rarr (C',u')$ are pairs of morphisms $f:C \rarr C'$ in the base $\cC$, and $g: u \mapsto u'$, being defined by $(Ff)(u) \rarr u'$ in $FC'$. It is this last definition we would like to discuss presently. We are looking at the following situation:

\beq
\setlength{\unitlength}{0.8cm}
\begin{picture}(5,6)
	\thicklines
	\put(-0.2,6){\line(1,0){0.4}}
	\put(-0.2,3){\line(1,0){0.4}}
	\put(0,6){\line(0,-1){3}}
	\put(4.8,6){\line(1,0){0.4}}
	\put(4.8,3){\line(1,0){0.4}}
	\put(5,6){\line(0,-1){3}}
	\put(-1,4.5){$FC$}
	\put(6,4.5){$FC'$}
	\put(0.1,4.9){$u$}
	\put(0.5,5){\vector(1,0){3.1}}
	\put(1.8,5.2){$Ff$}
	\put(3.8,4.9){$Ffu$}
	\put(4.5,4.8){\vector(0,-1){0.7}}
	\put(4.6,4.5){$g$}
	\put(4.3,3.7){$u'$}
	\multiput(0.5,4.8)(0.4,-0.1){9}{\circle*{0.1}}
	\put(4,3.9){\vector(3,-1){0.2}}
	\put(2.3,3.7){$g$}
	\multiput(0,2.8)(0,-0.3){5}{\circle*{0.07}}
	\put(0,1.7){\vector(0,-1){0.2}}
	\put(-0.2,1){$C$}
	\multiput(5,2.8)(0,-0.3){5}{\circle*{0.07}}
	\put(5,1.7){\vector(0,-1){0.2}}
	\put(4.8,1){$C'$}
	\put(0.5,1.2){\vector(1,0){4}}
	\put(2.3,0.5){$f$}
\end{picture} \nonumber
\eeq

Having the functor $Ff$, the morphism $u \rarr u'$ above is dictated by a choice of $g: Ffu \rarr u'$, so having one such morphism is all that we need. For that reason we will denote the morphism $u \rarr u'$ by the same letter $g$ by abuse of notation. Observe that one can define the manifestation of $Ffu$ in $FC'$ by:
\beq
\text{Man}_{FC'}(Ffu) = \Hom_{FC'}(Ffu,-) \nonumber
\eeq
and:
\beq
\Hom_{FC'}(Ffu,u') \ni g:Ffu \mapsto u' \nonumber
\eeq
There is an obvious forgetful functor $F_0: \int F \rarr \cC$. In some instances, one can reconstruct $F$ from the data of $\int F$ and $F_0$.

\section{Classical $\infty$-categorical Grothendieck Construction}

As pointed out above, one cannot always reconstruct the functor $F$ from its category of elements $\int F$, as it is sometimes called. The straightening/unstraightening construction  of \cite{Lu} is the $\infty$-categorical analogue of the traditional Grothendieck construction, as reviewed in Propositions 4.5 and 4.11 of \cite{A}, which further provides a setting where one can reconstruct the desired functors.\\

We start by first considering the statement:
\beq
\coCFib(\cC) \simeq \Fun(\cC,\Catinf) \nonumber
\eeq
from \cite{Lu}, where $\cC \in \Catinf$, as reviewed in \cite{AMG}. We aim to deconstruct this equivalence to the point where it appears we have a ground state from which one can develop two pictures, one being provided by the right hand side of the above equivalence, the other being provided by the left hand side.\\

\subsection{co-Cartesian fibrations}
We will start with $\coCFib(\cC)$ since this is the most interesting side of the equivalence we are considering from an information theoretic perspective. Recall that $\coCFib(\cC) = N((\SetplusD)_{/\cC})^{\circ}$ is the nerve of the subcategory of fibrant and cofibrant objects of the model category of marked simplicial sets over $\cC$ $(\SetplusD)_{/\cC}$ with the co-Cartesian model structure. We focus on the latter slice category of marked simplicial sets $(\SetplusD)_{/\cC}$. We work in full generality, and for that purpose we start by considering generic co-Cartesian fibrations $p:X \rarr S$ in $\SetD$. 

\subsubsection{Inner fibrations}
First and foremost, co-Cartesian fibrations are inner fibrations. By definition those have the RLP with respect to all inclusions $\Lambda^n_i \hrarr \Dn$ for $0 < i < n$, which one can represent as:
\beq
\xymatrix{
	\Lambda^n_i \ar@{^{(}->}[d]_{i} \ar[r]^{\sigma} &X \ar[d]^{p} \\
	\Dn \ar@{..>}[ur]^{\psi} \ar[r]_{\rho} & S
} \nonumber
\eeq
Inner horn inclusions $\Lambda^n_i \hrarr \Delta^n$ are regarded as probing $\Delta^n$ since in addition to being inclusions, the geometric realization $|\Lambda^n_i|$ is a retract of $|\Delta^n|$ in $\Top$. Moreover, owing to the definition of horns $\Lambda^n_i$ as being $\Delta^n$ from which we have removed its interior and its $i$-th face, one can characterize such probings as being peripheral. Having one solid diagram such as the one above means looking at $p$ from the perspective of an inner horn inclusion, thus we are focusing on the peripheral probing character of $p$. That we have the RLP at this point means $\sigma$ can be filled to $\psi$, something we further regard as providing a dynamic notion of density. Once maps are based at $\Dn$, we are looking at:
\beq
\xymatrix{
	&X \ar[dd]^p\\
	\Dn \ar[ur]^{\psi} \ar[dr]_{\rho} \\
	&S
} \nonumber
\eeq
To summarize, $p$ can be studied from the perspective of standard peripheral probings. From that viewpoint, $p$ being an inner fibration, one can then use the density of $\sigma$. Collecting things together, we are led to referring to $p$ as a \textbf{holographic probing map}. Simply said, $X$ can be regarded as fully probing $S$ from within.\\

Observe that if $p:X \rarr S$ is an inner fibration, the fibers $X_s = X \times_S \{s\}$ of $p$ over points $s \in S$ are $\infty$-categories, and morphisms between such fibers are correspondences, as explained in \cite{Lu}. This does not allow one to basculate to maps $S \rarr X$ yet. But to further ask that those maps $p$ be co-Cartesian fibrations will allow us to make such a move, as we will see below.\\

Next come the additional characteristic features of co-Cartesian fibrations, which require the notion of co-Cartesian edges, and those necessitate the introduction of a particular kind of simplicial set, which we discuss presently.\\

\subsubsection{Thickenings of maps}
Part of the definition of co-Cartesian edges makes use of a particular type of simplicial set defined in \cite{Lu} as follows. For any map $p:K \rarr S$ in $\SetD$, there is a simplicial set $S_{/p}$ characterized by the following universal property:
\beq
\Hom_{\SetD}(Y,S_{/p}) = \Hom_p(Y \star K, S) \nonumber
\eeq
where as defined in \cite{Lu}, the $p$ subscript on the right hand side Hom means we only consider those morphisms from the join $Y \star K$ that restrict to $p$ on $K$ alone. What this characterization means is that we enlarge our perspective on $p$ by considering all those maps from joins with $K$ that originate from $p$. If $Y \star K$ is an enlargement of the base, the collection of all maps with such a base viewed from the perspective of $Y$ provides $S_{/p}$, as given by the left hand side of the above equivalence. Consequently $S_{/p}$ can be viewed as a thickening of $p$, a way to actually turn a map of simplicial sets into a simplicial set, something that will prove useful when one deals with commutative diagrams in $\SetD$.\\

As a shorthand, if $\cC \in \Catinf$ and one considers $p: \Delta^0 \rarr \cC$ with image $X$, then as in \cite{Lu} one just denotes $\cC_{/p}$ by $\cC_{/X}$.\\

\subsubsection{co-Cartesian edges}
Those thickenings being defined one can now define co-Cartesian edges of a given map $p:X \rarr S$ in $\SetD$. If $p$ is an inner fibration, then an edge of $X$ is said to be $p$-co-Cartesian if it is $p$-Cartesian with respect to $p^{\op}: X^{\op} \rarr S^{\op}$, where $f:x \rarr y$ in $X$ is $p$-Cartesian if:
\beq
X_{/f} \rarr X_{/y} \times_{S_{/py}} S_{/pf} \nonumber
\eeq
is a trivial Kan fibration, in particular an equivalence of simplicial sets. One observation we can make at this point is that here we have a statement written from the perspective of $f$ in $X$. Further, in a homotopical context it makes sense to work with simplicial sets instead of just maps, hence the use of simplicial sets of the form $X_{/f}$. That we have such an equivalence above means that those edges that project to $py$ and $pf$ in the base ought to be $f$ up to equivalence, since $X_{/f}$ is a thickening of $f$ itself. In other terms being $p$-co-Cartesian ensures the uniqueness of the edge in $X$ being considered, up to a contractible ambiguity.\\

\subsubsection{co-Cartesian fibrations}
We finally arrive at the definition of co-Cartesian fibrations as defined in \cite{Lu}. $p: X \rarr S$ a co-Cartesian fibration is first an inner fibration, and as we argued after using a density argument this amounts to putting $X$ and $S$ on a same footing relative to $\Dn$, thereby presenting $X$ as it were as a holographic presentation of $S$. Then to say that $p$ is a co-Cartesian fibration means for any edge $f:x \rarr y$ in the base, for any $\tilde{x} \in X$ projecting to $x$ under $p$, meaning if we place ourselves at a point in $X$ from which we can have possible maps projecting down to $S$, then one is guaranteed to find a $p$-co-Cartesian edge $\tilde{f}: \tilde{x} \rarr \tilde{y}$ of $X$ projecting down to $f$ in the base, which really means we can lift $f$ to an edge of $X$ up to a contractible ambiguity, which again provides one with a sort of homotopical uniqueness. This allows one to decorate the maps of $S$ by those of $X$ and this in a unique manner. Equivalently said, asking that we work with co-Cartesian fibrations $p: X \rarr S$ therefore amounts to first probing $S$ by a dense substructure $X$ (which by the same token provides a holographic presentation of $S$), and morphisms in the base come from morphisms in such substructures, and this in a unique way. \\

Now if inner fibrations, with their correspondences as maps between fibers, did not allow one to fully define a functor from $S$ into $X$, co-Cartesian fibrations however allow one to do just that by virtue of the fact that lifts are up to a homotopical ambiguity. At this point one has completely deconstructed those maps of $(\SetplusD)_{/\cC}$ that contribute to $\coCFib(\cC)$ and we find ourselves in a situation where we can either reconstruct what we have back to $(\SetplusD)_{/\cC}$, or basculate to a picture where one considers functors from $\cC$ into $\Catinf$ instead, in accordance with the adjunction equivalence $\coCFib(\cC) \simeq \Fun(\cC,\Catinf)$ of \cite{Lu}.\\

\subsection{Functors into $\Catinf$}
Given $\cC \in \Catinf$, $\Fun(\cC,\Catinf)$ provides all possible manifestations of $\cC$ into $\Catinf$, but most noteworthy, $F \in \Fun(\cC,\Catinf)$ being one such functor it attaches to each point $C \in \cC$ an $\infty$-category $F(C) \in \Catinf$, and this functorially, so effectively we attach a microstructure to $\cC$, and this in a coherent fashion. Collectively, $\Fun(\cC,\Catinf)$ therefore provides all the possible ways to do just that.\\

Technically, $\Fun(\cC,\Catinf) = \Fun(\cC,N((\SetplusD)^{\circ})) = N((\SetplusD)^{\gC[\cC]})^{\circ})$ by Proposition 4.2.4.4 of \cite{Lu}.

\subsection{Equivalence}

One starts with the $\infty$-categorical Grothendieck construction, Theorem 3.2.0.1 of \cite{Lu}, which states that if $K \in \SetD$, $\cC \in \CatD$, $\phi: \gCK \rarr \cCop$ is an equivalence in $\CatD$, then one has a Quillen equivalence:
\beq
\SetDplusK \adj^{\Stphiplus}_{\Unphiplus} \SetDplusC \label{QuillEquiv}
\eeq
where on $\SetDplusK$ one puts the cartesian model structure, and one puts the projective model structure on $\SetDplusC$. It is not in this form though that an $\infty$-categorical Grothendieck construction is self-evident. To make it apparent, one has to take the underlying $\infty$-categories of the simplicial categories in \eqref{QuillEquiv} as summarized in \cite{AMG}. While doing so, one invokes Proposition 4.2.4.4 of \cite{Lu} which provides a categorical equivalence of simplicial sets $N((\cU^{\cC})^{\circ}) \xrarr{\simeq} \Fun(K, N(\cU^{\circ}))$ if $\gCK \xrarr{\simeq} \cC$. One uses this with $\cU = \SetplusD$, along with the equivalence $(\SetDplus)^{\circ} \simeq \Catinf^{\Delta}$, $\Catinf^{\Delta}$ the simplicial category of small $\infty$-categories, whose simplicial nerve $N(\Catinf^{\Delta}) = \Catinf$, $\Catinf$ the $\infty$-category of small $\infty$-categories. This gives us $N((\SetDplus)^{\cC})^{\circ} \simeq \Fun(K, \Catinf)$. But the Quillen equivalence \eqref{QuillEquiv} also provides an adjunction equivalence $N(\SetDplusK)^{\circ} \adj N(\SetDplusC)^{\circ} \simeq \Fun(K, \Catinf)$ as a consequence of Lemma 3.2.4.1 of \cite{Lu}. Note that on $(\SetplusD)_{/K}$ one now has the co-Cartesian model structure. That we can do this is warranted by taking opposites of categories, once on $(\SetplusD)_{/K}$, and once on $(\SetplusD)^{\cC}$, as reviewed in \cite{A}. Thus one arrives at the (dual) equivalence of $\infty$-categories as recounted in \cite{AMG} for $K = \cC \in \Catinf$: $\coCFib(\cC) = N((\SetplusD)_{/\cC})^{\circ} \xrarr{\simeq} \Fun(\cC,\Catinf)$, the former $\infty$-category being the $\infty$-category of co-Cartesian fibrations on $\cC$, since on $(\SetplusD)_{/\cC}$ one has put the co-Cartesian model structure, for which fibrant objects are co-Cartesian fibrations.\\

Phenomenologically, this equivalence $\coCFib(\cC) \simeq \Fun(\cC, \Catinf)$ has the following interpretation. As we argued, to consider co-Cartesian fibrations over $\cC$ means looking at maps $p: X \rarr \cC$ for which whenever we have edges in the base $f: x \rarr y$ that could come from edges in $X$ in the sense that there is some $\tilde{x}$ for which $p\tilde{x} = x$, then there is indeed such an edge $\tilde{x} \rarr \tilde{y}$ in $X$ mapping down to $\cC$, and this in a unique way. In accordance with our interpretation of co-Cartesian fibrations, edges of $\cC$ come from within by density, they are the result of probing the structure of $\cC$. More specifically, they arise as holographic presentations of edges of $X$. That the $\infty$-category of co-Cartesian fibrations is equivalent to $\Fun(\cC,\Catinf)$ says that probing the structure of $\cC$ in this fashion produces an $\infty$-categorical structure at each object of $\cC$, something we regard as putting a microstructure on the objects of $\cC$, and to map those microstructures together by functoriality. $\Fun(\cC,\Catinf)$ provides all the ways to do just that. Simply said, holographically probing $\cC$ results in revealing its $\infty$-categorical microstructure.\\

What we have achieved is that by progressively clarifying the structure of $\coCFib(\cC)$, we eventually reached a stage where one can either reconstruct what we have back into $\coCFib(\cC)$, or basculate to another picture $\Fun(\cC,\Catinf)$. This equilibrium state, a \textbf{ground state} so to speak, encapsulates both macroscopic points of view as potentialities from an information theoretic perspective. We can represent this diagramatically as follows:
\beq
\xymatrix{
	\coCFib(\cC) \ar[dr]_{dec} \ar@{.>}[rr]^{\simeq} && \Fun(\cC,\Catinf) \\
	& \text{Ground} \ar[ur]_{rec}
} \nonumber
\eeq
with an equivalent diagram with all arrows reversed. One map is an information theoretic \textbf{deconstruction}, we have fully dissected $\coCFib(\cC)$. Once this is achieved we have reached a ground state where two dual, macroscopic pictures coexist. A possible \textbf{reconstruction} provides $\Fun(\cC,\Catinf)$. From this perspective perceptions and manifestations appear as two dual realizations of a same ground state. Relative to $\cC$, we have a conservation law: the algebraic information that comes into it via co-Cartesian fibrations is also given by the microstructure it generates via functors $\cC \rarr \Catinf$.\\

It is worth pointing out that the ground state is not something that one can define explicitly in the same manner that $\coCFib(\cC)$ was defined for example. Rather, the process of straightening/unstraightening the model categories from which either $\infty$-category originates unravels their structure to the point where we can basculate from one model category to the other. At that moment the ground state has been made manifest. \\

\section{Grothendieck construction for partial functors}
It is in \cite{A} that the generalization of the $\infty$-categorical construction to partial functors is considered. By partial functor we mean functors $\cD \rarr \Catinf$ that are not defined on all morphisms of $\cD$. It is interesting that this lack of information defines a particular class of object by complementarity. In other terms if $\Fun(\cD,\Catinf)$ defines an $\infty$-category, $\Fun_{\text{part}}(\cD,\Catinf)$ defines yet another $\infty$-category. The proper way to discuss partial functors is by using categorical patterns, which we define first, followed by fibered categories in patterns, which provide a complete and efficient way to characterize partial functors (see \cite{Lu2}). One can come up with a Grothendieck construction-like categorical equivalence between two $\infty$-categories as done in \cite{A}, a generalization of the straightening/unstraightening construction of \cite{Lu} to the case of partial functors $\cD \rarr \Catinf$.

\subsection{Categorical patterns}
We follow \cite{A}, which itself is based on \cite{Lu2}. We do not need the full definition of categorical pattern but only that of a commutative categorical pattern, which \cite{A} just refers to as categorical pattern, and we will follow his lead. Let $\cD \in \Catinf$ being fixed. By Definition 2.14 of \cite{A}, a categorical pattern $\bp$ on $\cD$ is a pair $\bp = (\ED, \{\pa: \KaT \rarr \cD\}_{\alpha \in A})$ consisting of a set of edges $\ED$ of $\cD$ that contains all equivalences and degenerate edges, making $\overline{\cD} = (\cD, \ED)$ into a marked simplicial set, and maps $\pa: \KaT \rarr \cD$ that are such that $\pa(\KaT) \subset \ED$.

\subsection{Fibered categories}
As pointed out in \cite{A}, partial functors $(\cD,\bp) \rarr \Catinf$ can be conveniently transcribed into fibrations over $(\cD,\bp)$, by using what is called $\bp$-fibered categories over $\cD$.\\

One works with $\SetplusD/_{\overline{\cD}}$, simply denoted by $\SetplusDp$ as in \cite{A}, and which one refers to as the category of marked simplicial sets over $\bp$. We will only consider those objects of $\SetplusDp$ that allow one to basculate to partial functors by way of an equivalence, in the same manner that we had the equivalence $\coCFib(\cD) \simeq \Fun(\cD,\Catinf)$. If this equivalence is any indication we want objects $\overline{X} = (X,\EX) \rarr \overline{\cD} = (\cD,\ED)$ of $\SetplusDp$ to have the flavor of co-Cartesian fibrations. One can build on the work presented in the previous section and for that reason one asks that the underlying maps of simplicial sets $X \rarr \cD$ be inner fibrations. This is a prerequisite, and this allows us to eventually have functors from $\cD$ valued in $\infty$-categories. Recall that having inner fibrations means holographically probing the structure of $\cD$, and the result of such an operation is to attach $\infty$-categories at each point of $\cD$.\\

Beyond this point one has to impose conditions that take into consideration the fact that our simplicial sets are decorated by sets of edges, and that in addition those are imbued with a pre-dynamics in the base by way of the maps $\pa$. We will cover each of these conditions in separate sections, and together they define what \cite{A} refers to as $\bp$-fibered objects of $\SetplusDp$.\\

\subsubsection{co-Cartesian pullbacks}
Let $\overline{p}: \overline{X} \rarr \overline{\cD}$ be a marked simplicial set of $\SetplusDp$ with underlying map $p:X \rarr \cD$. As a first condition for $\overline{p}$ to define what \cite{A} refers to as a $\bp$-fibered simplicial set, one asks that $p$ be an inner fibration as pointed out above. The second condition of Definition 2.1 where $\bp$-fibered simplicial sets are defined now asks that for any edge $\tilde{e} \in \EX$, the induced map $p': X \times_{\cD} \Delta^1 \rarr \Delta^1$ be a co-Cartesian fibration. This should be understood in the following sense:
\beq
\xymatrix{
	X \times_{\cD} \Delta^1 \ar[d] \ar[r]^-{p'} & \Delta^1 \ar[d]^{e} \ar@{.>}[dl]_{\tilde{e}} \\
	X \ar[r]_p & \cD
}\nonumber
\eeq
the edge $e = p \tilde{e}$ of $\cD$ has a holographic presentation in $X \times_{\cD}\Delta^1$ via $p'$.\\

\subsubsection{locally $p$-co-Cartesian edges}
Since $\overline{X} = (X,\EX)$ and $\overline{D} = (\cD, \ED)$, a map of marked simplicial sets $\overline{X} \rarr \overline{\cD}$ maps $\EX$ into $\ED$, thus condition (3) of \cite{A} for the definition of $\bp$-fibered marked simplicial sets is a given, if $\tilde{e} \in \EX$, it ought to map to $\ED$. That the converse is true allows one to lift marked edges of $\cD$ up to $X$, and for this to be well-defined we want $\tilde{e}$ to be locally $p$-co-Cartesian, that is it is $p'$-co-Cartesian where $p'$ was introduced in the previous subsection. What that really means is focusing on $e = p\tilde{e}$, $\tilde{e}$ is the only lift possible at the level of $X$ up to a contractible ambiguity.\\

The other conditions defining a $\bp$-fibered marked simplicial set over $\cD$ are less transparent. For the sake of interpreting the $\infty$-categorical Grothendieck construction in the present setting they can be ommitted, though one has to keep in mind that they are needed in more technical treatments of this material. If anything, one can see those conditions as allowing us to treat partial functors as $\bp$-fibered marked simplicial sets over $\cD$.\\

\subsubsection{$\infty$-category of $\bp$-fibered simplicial sets}
Theorem 2.10 of \cite{A} puts a model structure on $\SetplusDp$ for which fibrant objects are $\bp$-fibered objects, which justifies calling the underlying $\infty$-category $N(\SetplusDp)^{\circ}$ the $\infty$-category of $\bp$-fibered objects, denoted $\pFib$. Here $N$ is the homotopy coherent nerve functor, right adjoint to the functor $\widetilde{\gC}$, conjugate to the functor $\gC$ used in \cite{Lu}.\\

\subsection{Bundles}
\subsubsection{Motivational definition of $\bp$-bundles}
Arakawa looked at functors from some $S \in \SetD$ into $\Fun_{\text{part}}(\cD,\Cat)$ in its $\pFib \in \Catinf$ incarnation, that is we consider $\Fun(S, \pFib)$, and he showed that one can transcribe that Grothendieck construction-style into something of the form $\cE_{/S}$ for some object $\cE$ to be defined. Let $\cD \in \Catinf$ fixed, $\bp = (\ED,\{\pa:\KaT \rarr \cD\}_{\alpha \in A})$ a fixed categorical pattern on $\cD$, resulting in $\pFib$. We anticipate to have $\cE_{/S}$ being a bundle over $S$ with fibers in $\pFib$. To force the issue, one considers $\pFib$-fibered objects over $S$. One can start by a simple bundle over $S$. If $\cX \in \cE$, and we tentatively write $p:\cX \rarr S$ for an object of $\cE_{/S}$, then one asks that for all $s \in S$, $\cX_s \in \pFib$, and for any $f:s \rarr s'$ in $S$, one has an induced morphism $f_{!}: \cX_s \rarr \cX_{s'}$ of $\bp$-fibered objects. \\

Next, by virtue of the fact that $\pFib = N(\SetplusDp)^{\circ}$, and given that fibrant objects for the model category $\SetplusDp$ are $\bp$-fibered objects $X \rarr \cD$, one has $\cX_s$ of the form  $\{X \rarr \cD\}$ for one such object, and this for all $s \in S$. To have a map $\cX \rarr S$, one starts by considering a diagram of the form:
\beq
\xymatrix{
	X \ar[dr]_{q} \ar[rr]^p && S \times \cD \ar[dl]^{\pi} \\
	& S
} \nonumber
\eeq
and the bundle conditions we just gave now read $X_s = X \times_S \{s\} \rarr \cD$ is $\bp$-fibered, and $X_s \rarr X_{s'}$ is a morphism of $\bp$-fibered objects. For the sake of basculating from bundles of this sort to having functors into $\pFib$, one asks that $q$ be a co-Cartesian fibration. One also wants a thickened version thereof for $p:X \rarr S \times \cD$. It is natural to at least ask that $p$ be an inner fibration, and one also asks that $p$ be co-Cartesian modulo $\cD$, which means for all $x \in X$ with $px = (s,D) \in S \times \cD$, for any $(f,g): (s,D) \rarr (s',D')$, if $g$ is an equivalence in $\cD$, then there is a $p$-co-Cartesian edge $\tilde{e}:x \rarr x'$ with $p\tilde{e} = (f,g)$. This condition is needed as $\pi \circ p = q$. If all those conditions are satisfied, one has a $\bp$-bundle over $S$ as defined in Definition 3.1 of \cite{A}.\\

\subsubsection{$\infty$-category of $\bp$-bundles over $S$}
As usual since \cite{Lu}, one adopts the notation $S$ to denote the marked simplicial set $S^{\sharp}$, and in particular when one deals with $S^{\sharp} \times \bp$-fibered objects, those will simply be referred to as $S \times \bp$-fibered objects. By Proposition 3.5 of \cite{A}, $\bp$-bundles $p:X \rarr S \times \cD$ as defined above are the $S \times \bp$-fibered objects of $\SetplusDSp$. By Corollary 3.6 of \cite{A} the model structure we put on the latter category is the one for which fibrant objects are the $\bp$-bundles over $S$ whose marked edges are the locally $\underline{p}$-co-Cartesian edges of $\underline{p}:X \rarr \cD$. This justifies the definition:
\beq
N(\SetplusDSp)^{\circ} = \pBun(S) \nonumber
\eeq

\subsection{Equivalence}
\subsubsection{Organizational commutative diagram}
We argue basculating in accordance to the main result $\pBun(S) \simeq \Fun(S,\pFib)$ of \cite{A} hinges on hitting a ground state in the classical $\infty$-categorical Grothendieck construction of \cite{Lu}, and on considering the corresponding induced maps on $\bp$-fibered or $\bp$-bundle objects, at which point we preserve a notion of ground state. From there one can go from bundle representations to manifestations into $\pFib$ in accordance with the equivalence mentioned above. The other equivalences involved in showing such a result do not use a duality in the same manner that the straightening/unstraightening construction did. We are looking at the following diagram:
\beq
\xymatrix{
	(\SetplusD)^{\tilde{\gC}[S]} \ar@{~>}[d]  \ar[r]^{ \widetilde{\text{Un}}^+} & \SetplusDS \ar@{~>}[d] \\
\SetplusDptCS \ar[d]_{\infty} \ar[r]^{\simeq}_{\text{Thm 5.7}} & \SetplusDSUnp \ar[d]_{\infty} \ar[r]^{\simeq}_{\text{Prop 5.5}} & \SetplusDSp \ar[d]_{\infty} \\
N(\SetplusDptCS)^{\circ} \ar[d]_{\simeq} \ar[r] & N(\SetplusDSUnp)^{\circ} \ar[r] &N(\SetplusDSp)^{\circ} \ar[d]|-{=}\\
\Fun(S,N(\SetplusDp)^{\circ}) \ar[d]|-{=} && \pBun(S)\\
\Fun(S,\pFib) \ar@{.>}[urr]_{\simeq}
} \nonumber
\eeq

The top equivalence $\widetilde{\text{Un}}^+$ we have by virtue of \cite{Lu}, and this is really where we have the ground state we mentioned in the previous section appearing, where both pictures coexist. The crux of \cite{A} is that this is preserved next line at the level of $\bp$-fibered marked simplicial sets over $\bp$ by virtue of Proposition 2.16, and an astute use of a combination of the less obvious Proposition 3.5 and Lemma 2.23, along with the properties of Quillen functors, all with an aim to not directly involve $\bp$-fibrations. The first equivalence on the second line is given in Theorem 5.7. It uses the equivalence $\cD \xrarr{\simeq} \widetilde{\text{Un}}_{\Delta^0}(\cD)$, which on the categorical pattern $\bp$ produces what we denote by $\Unp = (\text{Un} \ED, \{ \text{Un}\pa\}_{\alpha \in A})$. The equivalence that follows on the same line is given by Proposition 5.5, and both are Quillen equivalences. The main point is that the ground state which is present in the classical $\infty$-categorical Grothendieck construction generalizes to the case of dealing with categorical patterns, albeit in a less transparent manner.\\

To come back to Theorem 5.7, to argue in favor of the existence of the functor $\SetplusDptCS \rarr \SetplusDSUnp$, it suffices to start with simplicial functors valued in $\SetplusDp$. This means considering those functors that are really of the form $F \rarr \delta{\overline{\cD}}$, where $\delta(\overline{\cD})$ is the constant $\overline{\cD}$-valued functor, with associated unstraightening that morally looks like:
\beq
	\Un F  \rarr \Un \delta(\overline{\cD}) = S^{\sharp} \times \Un (\overline{\cD})  \nonumber
\eeq
object of the category $\SetplusDSUnp$, hence the functor $\SetplusDptCS \rarr \SetplusDSUnp$ is well-defined.\\

\subsubsection{Ground state for partial functors}
As a first motivation for the existence of a ground state in the present picture whereby $\pBun(S) \simeq \Fun(S, \pFib)$, consider a $\bp$-bundle over $S$:
\beq
\xymatrix{
	X \ar[dr]_q \ar[rr]^p &&S \times \cD \ar[dl]^{\pi}\\
	&S
} \nonumber
\eeq
where $q$ is a co-Cartesian fibration, in particular an inner fibration, and so is $p$. It turns out the fibers of this bundle are in $\pFib$, thus we aim to show we have a transition of the form:
\beq
\xymatrix{
	\pFib \ar[d] \\
	s \in S &\ar@{~>}[r] & s \mapsto \psi(s) \in \Ob(\pFib)
} \nonumber
\eeq
for some functor $\psi$. That we preserve fibers will point to the existence of a ground state. To argue in favor of this transition, it suffices to show that the straightening map $\tSt$ maps $\pFib \rarr s$ to $s \mapsto \Ob(\pFib)$. In \cite{A}, Arakawa shows that we have an equality $\tStplusphi = \tStpplusphi$ of more general straightening functors, with $\phi: \tgCS \rarr \cC$ a given simplicial functor, and where $\tStpplusphi$ is Lurie's version of the straightening functor. Recall that if $p: \overline{X} \rarr S^{\sharp}$ is a map of marked simplicial sets, then one defines $\tStphi'(p)(C) = (\tgC[\Xlhd] \coprod_{\tgC[X]} \cC)(\infty,C) \in \SetD$, where $\infty$ is the cone point of $\Xlhd$. Diagrammatically:
\beq
\xymatrix{
	X \ar[d]_p \ar[r] & \tgC[X] \ar[d] \ar[r] & \tgC[\Xlhd] \ar[dd] \\
	S \ar[r] &\tgCS \ar[d]_{\phi} \\
	&\cC \ar[r] &\cC\coprod_{\tgC[X]} \tgC[\Xlhd]
} \nonumber
\eeq
Thus $\tStphi'(p)(C)$ gives the simplicial set of morphisms from $\infty$ in the fiber to $C$. For the sake of having an equivalence at the level of $\infty$-categories, we take $\phi$ to be an equivalence in $\CatD$. Then $C \in \Ob(\cC) \simeq \Ob(\tgCS) = S_0 \ni s$, so we do have something of the form $X_s \rarr s$ mapping to $\phi(s) \rarr \SetD$. Dealing with markings now we need a simplicial functor $\cC \rarr \SetplusD$ and for that purpose one takes the minimal marking on $\tStphi(p)(C) = \tStphi'(p)(C)$ for which $\tStplusphi(p)$ achieves that, with the added condition that any $u:x \rarr y$ in $\overline{X}$ marked maps to $\tilde{u}$ in $\tStphi(p)(\phi py)$, marked as well, in accordance with the following set up:
\beq
\xymatrix{
	\tgC[(\Delta^1)^{\lhd}] \ar[r] & \tgC[\Xlhd] \ar[r] & \tgC[\Xlhd] \coprod_{\tgC[X]} \cC \ni \tStphi(p)(\phi py) \\
	u:\Delta^1 \ar[u] \ar[r] \ar@{.>}[urr] & X \ar[u]
} \nonumber
\eeq
in such a manner that we preserve marked simplicial sets, as presented in \cite{A}. Finally we argue we preserve objects of $\pFib$, which was not obvious so far. This is made manifest once we take the nerve of the subcategories of fibrant and cofibrant objects, the reason being that $\pFib$ objects are preserved from the perspective of the straightening functor, not by construction, but by adjunction with the unstraightening functor which does preserve $\pFib$-like objects.\\

That we have such a transition whereby we preserve some structure such as fibered objects is proved in general in \cite{A} by considering the unstraightening functor instead, where $\phi: \tgCS \rarr \cC$ is a morphism of simplicial categories (though for our purposes $\phi = id$ is sufficient), $\bp$ is a commutative categorical pattern on some $\cD \in \Catinf$, $\tUnpphi: (\SetplusDp)^{\cC}_{\text{proj}} \rarr \SetplusDSUnp$, a right Quillen functor by Proposition 5.6 of \cite{A}. Before putting model structures on the underlying categories of interest, we find ourselves in the same situation as with the straightening functor, we have simplicial sets mapping to simplicial sets but $\pFib$-like objects are not apparent yet. For that to happen one has to add the model structures in the picture and in particular one has to put the projective model structure on $(\SetplusDp)^{\cC}$. In the proof of Proposition 5.6 one shows that it is sufficient to just consider $S \in \Catinf$ and $\phi = id$. Then $\tUnpphi$ preserves fibrant objects. In the projective model structure on $\SetplusDptCS$, an object $F$ therein is fibrant if $F(C)$ is fibrant in $\SetplusDp$ for all $C \in \tgCS$, and those are the objects of $\pFib$, hence $F \in \Fun(S,\pFib)$. $F$ is mapped by $\tUnpphi$ to a fibrant object of $\SetplusDSUnp$, which means in particular by Proposition 3.5 that it is a $\bp$-bundle over $S$, in particular it has $\pFib$-valued fibers.\\

To focus on the preservation of $\pFib$-objects is the analogue in the present situation of preserving the $\infty$-categorical fibers of $\coCFib(\cC)$ while making the transition to $\Fun(\cC,\Catinf)$, the dynamics of which highlights the presence of a ground state. Here we have an added layer of complexity because we also have to track $\bp$-like data.\\

It appears a deconstruction of $\pBun(S)$ and a reconstruction into $\Fun(S,\pFib)$ hinges on the presence of a ground state in the traditional $\infty$-categorical Grothendieck construction of Lurie, in conjuction with the use of properties of Quillen functors and a variety of results proved in \cite{A}, in particular Proposition 3.5, Proposition 2.16 and lemma 2.23. Collecting things together, algebraically we have:
\beq
\xymatrix{
	\pBun(S) \ar[dr]_{dec} \ar@{.>}[rr]^{\simeq} && \Fun(S,\pFib) \\
	& \text{Ground} \ar[ur]_{rec}
} \nonumber
\eeq
with a similar triangle with arrows reversed.\\

\section{Phenomenology}
\subsection{Ground state formalism}
Whether it be in the traditional $\infty$-categorical Grothendieck construction or the one used for partial functors, one can illustrate the presence of a ground state by focusing on the preservation of fibers as pointed out above. In the former case we have an equivalence $\coCFib(\cC) \simeq \Fun(\cC, \Catinf)$. If $p:X \rarr \cC$ is a co-Cartesian fibration, in particular an inner fibration, its fibers are $\infty$-categories. Let $F$ be its associated functor from $\cC$ into $\Catinf$. We have:
\beq
\xymatrix{
	X_C \in \Ob(\Catinf) \ar[d]_{p} \\
	C \in \Ob(\cC) \ar@{~>}[r] & C \mapsto F(C) \in \Ob(\Catinf) 
} \nonumber
\eeq
For partial functors we have an analogous picture. It turns out the fibers of objects of $\pBun(S)$ are valued in $\pFib$, thus we are looking at:
\beq
\xymatrix{
	\pFib \ar[d] \\
	s \in S \ar@{~>}[r] & s \mapsto \psi( s) \in \Ob(\pFib)
} \nonumber
\eeq
for some functor $\psi: S \rarr \pFib$.\\

In both cases therefore we have the $\infty$-category of $\cE$-valued presheaves from $S^{\op} \in \SetD$ being equivalent to an $\infty$-categorical bundle over $S$ with fibers in $\cE$, the bundle being $\coCFib(S)$ (with $\cE = \Catinf$) in one case and $\pBun(S)$ (with $\cE = \pFib$) in the other. A close examination of the dynamics at the level of fibers over points of $S$ and their associated maps from $s \in \Ob(S)$ on the other side of the equivalence makes manifest the existence of a ground state for such equivalences. Formally:\beq
\xymatrix{
	\cE \text{-Bun}(S) \ar[dr]^{dec} \ar@{.>}[rr]^{\simeq} && \Fun(S,\cE) \\
	& \text{Ground} \ar[ur]_{rec}
}\nonumber
\eeq
with a similar diagram with all arrows reversed. A suggestive presentation as:
\beq
\cE \text{-Bun}(S) \rlarr \text{Ground} \rlarr \Fun(S,\cE) \nonumber 
\eeq
displays the ground state as this intermediate $\infty$-category, a sort of middle ground, where the informational content of both theories is present. Moving in either direction away from this ground state breaks its symmetry and gives one picture or the other.\\

We can abstract this further as follows: starting with two formalisms $A$ and $B$ that are dual to each other, with an equivalence $E_{AB}$ between them, we have:
\beq
\xymatrix{
	A \ar@{-}[dr]_{rec/dec} \ar@{-}[rr]^{E_{AB}} && B \ar@{-}[dl]^{dec/rec} \\ 
	& \text{Ground}
} \nonumber
\eeq
where the ground state contains all the information pertaining to $A$ and $B$ and the equivalence between them. As we have seen in the case of the Grothendieck construction, we can always deconstruct either $\infty$-category $\coCFib(\cC)$ or $\Fun(\cC,\Catinf)$ into its fundamental constituents. This provides us with the deconstruction maps above. The reconstruction maps are the reverse processes. The equivalences between both $\infty$-categories partly arise from deconstructing either theory. The ground state also contains such information, so $E_{AB}$ itself appears therein in deconstructed form, though the above diagram does not display such a deconstruction.\\

The two theories $A$ and $B$ being equivalent, there is a core content $C_{AB}$ common to both, which is part of the ground state. Using this notation, focusing solely on the core content $C_{AB}$ above which we have either theory, the above diagram reads:
\beq
\xymatrix{
	A \ar@{-}[dr]_{rec/dec} \ar@{-}[rr]^{E_{AB}} && B \ar@{-}[dl]^{dec/rec} \\ 
	& C_{AB}
} \nonumber
\eeq
Denote by $N_A$ the nature of $A$, which encodes the information defining the algebraic structure of A. $N_A$ is needed to define the reconstruction map from $C_{AB}$ into $A$. Since $N_A$ is pertaining to $A$ only, $E_{AB}$ is some additional information, since it also takes account of the nature $N_B$ of $B$. Altogether the above triangle can be compactly represented as $(C_{AB}, N_A, N_B, E_{AB})$. All this information is part of the ground state, at which level any theory has dissolved into its constituting elements.\\

With regards to $E_{AB}$ containing additional information, for illustrative purposes consider the case of the $\infty$-categorical Grothendieck construction of Lurie, where the adjoint equivalence between $A = \coCFib(\cC)$ and $B = \Fun(\cC,\Catinf)$ arises from a Quillen equivalence of model categories. If the model categories capture part of the natures $N_A$ and $N_B$ of $A$ and $B$ respectively, and the deconstruction and reconstruction maps provide us with maps between the said model categories, the Quillen property that comes into the definition of $E_{AB}$ however brings in additional information, strictly speaking it is one property of maps of the ambient 2-category of model categories. Thus the above diagram is not meant to be understood as a commutative diagram, it is not sufficient to just deconstruct the theories to have the equivalences appear, more work may be needed in general.\\

To come back to the ground state itself, and if we use the $\infty$-categorical Grothendieck construction as impetus, we obtain a ground state in which the information content of either theory is in deconstructed form $(C_{AB},N_A,N_B,E_{AB})$, but not expressed in its most fundamental bits of information. In other terms what we are dealing with is not a fundamental ground state but an intermediate one. Thus once dual theories are captured together by a common ground, one can still refine that ground state to obtain something more fundamental, something we can picture as follows:
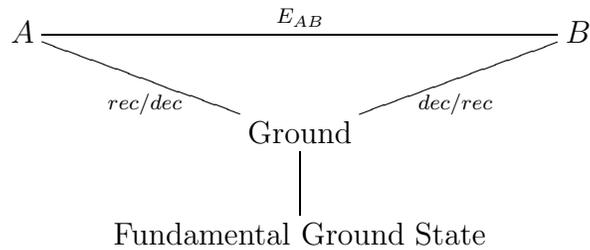
\beq
\xymatrix{
	A \ar@{-}[dr]_{rec/dec} \ar@{-}[rr]^{E_{AB}} && B \ar@{-}[dl]^{dec/rec} \\ 
	& \text{Ground} \ar@{-}[d] \\
	& \text{Fundamental Ground State}
} \nonumber
\eeq

Going from some intermediate ground to one more fundamental amounts to further deconstructing the information we have into more elementary chunks. For instance say at some intermediate stage we consider Quillen functors between model categories. At the next stage below one would have to deconstruct this: say we consider a left Quillen functor, in particular it preserves cofibrations among other things. One would have to explain at that level what it means to be a left Quillen functor, in particular, what ``preserving" means, and what are cofibrations. At the next level down from there, one would have to explain that to preserve something is a property, and that cofibrations is a particular type of map in a model category. At the next level, one would have to explain what is a property for a map, and what makes maps cofibrations. Next, one would have to define what is a map, and what is the relation between properties and maps. Then one has to explain what is the nature of a relation, etc...\\

An immediate question is why bother finding such a ground state? If one aims to come up with an all encompassing theory that captures all phenomena, it is important to come up with one formalism within which all concepts are put on a same footing, hence are undifferentiated, and this is implemented by dissecting the information we have into their most fundamental bits. In contrast, as soon as specific notions appear, distinctions can be made, and consequently one can argue we are no longer dealing with the ultimate ground state from which all theories can be derived, hence the importance of a search for a fundamental ground state where one is dealing with undifferentiated bits of information.\\

Such a ground state formalism has important physical implications as it points to the existence of a vacuum from which various dual theories can arise by \textbf{symmetry breaking}, which mathematically consists in condensing the information content of the ground state into either of the dual theories under consideration. We will come back to this point in a forthcoming work.

\subsection{Formal uncertainty principle}
Let us start with a traditional uncertainty principle $\Delta a \Delta b \geq \hbar$, $a$ and $b$ being issued from two respective theories $A$ and $B$. To write such an inequality presupposes that both $A$ and $B$ are defined in full, and it indicates that one is looking at an admixture of $a$ and $b$, the differential $\Delta a$ for instance being a measure of the amount of attention that is given to $a$. We can generalize that to theories as a whole as follows.\\

In the above ground state formalism, relative to $C_{AB}$, since $A$ and $B$ are dual theories, with common content encoded in $C_{AB}$, the more information from $C_{AB}$ is used to develop $A$, the less of the necessary information is left to develop $B$. Denote by $\Delta A$ the amount of information left to define $A$ in full relative to what has been constructed so far from $C_{AB}$. Denote this work in progress by $A_*$. One can formally write $\Delta A = A - A_*$. One can regard this as the amount of focus that is granted to $A$. As just argued, if $A$ and $B$ are two dual theories, the smaller is $\Delta A$, that is the closer we are to completing the construction of $A$, the less informational content is left in $C_{AB}$ to construct $B$, making $B_*$ small, so $\Delta B = B - B_*$ will be large. We can represent this by a \textbf{formal uncertainty principle} $\Delta A \Delta B \geq l > 0$, where $l$ denotes a lower bound, indicating that the informational contents of $A$ and $B$ vary in sync, which in itself is a corollary of having a duality. Note that these inequalities are formal and we are not dealing with quantifiers at the moment. Such inequality is just meant to indicate that there is a non-trivial lower bound to the joint information differentials between $A$ and $B$. We can also observe that in contrast to a traditional uncertainty principle that uses quantifiers, what we are dealing with here is discrete in nature since constructions/deconstructions are performed gradually. Thus we can regard formal uncertainty principles as being fractional in nature, as opposed to being continuous the way traditional uncertainty principles are. In the event that one has a bona-fide uncertainty principle with quantifiers, this obviously implies a formal uncertainty principle, the reason being that if one has a traditional uncertainty principle regarding two theories $A$ and $B$, those have to be well-defined. If that is the case they can be deconstructed, which provides us with a formal uncertainty principle. \\

Conversely one may ask to what extent does a formal uncertainty principle imply a regular one. We will conjecture below that a formal uncertainty principle implies having a duality. However a duality does not by any means imply having a traditional uncertainty principle as we now argue. Indeed, starting from dual theories $A$ and $B$, one can deconstruct both to some common content $C_{AB}$, and if there is deconstruction, there is also reconstruction. This provides us with the nabla-style diagrams. The composites of deconstruction and reconstruction maps are progressive and show that we have a formal uncertainty principle that comes for free. The horizontal edge corresponds to the initial equivalence between both theories, which provides an instantaneous change from $A$ to $B$ and vice-versa. There is nothing gradual about this change. Even if one deconstructs $E_{AB}$, to talk about a formal uncertainty principle at the level of equivalences makes no sense since it involves both $A$ and $B$. One could consider a continuous change illustrated by a traditional uncertainty principle involving elements of $A$ and $B$, thereby providing us with admixtures of $A$ and $B$-type objects at any stage of the transition from $A$ to $B$ and back in accordance with $E_{AB}$. However such an uncertainty principle between $A$ and $B$ is not a consequence of having a nabla-style diagram resulting from having $A$ and $B$ being dual. In other terms duality does not necessarily imply having a traditional uncertainty principle. However we at least suppose:
\begin{conj}
	If two theories $A$ and $B$ satisfy a formal uncertainty principle $\Delta A \Delta B \geq l > 0$, they correspond to dual theories.
\end{conj}

Starting from a formal uncertainty principle, the conjecture indicates that if say $A$ is developed, $B$ on the other hand is further away from being clearly developed, pointing to the existence of a joint information content $C_{AB}$, at which point both theories coexist in equal proportion. Thus one can write:
\beq
A \rlarr C_{AB} \rlarr B \nonumber
\eeq
with both dual arrows being dec/rec maps, hence invertible. By composition this gives us invertible maps from $A$ to $B$, hence the desired duality. Thus the conjecture is more a statement about dec/rec maps. Once both $A$ and $B$ are fully deconstructed, it is not automatic they share a common content $C_{AB}$. This is so if their individual deconstructions have subparts $C_A$ and $C_B$ from which both $A$ and $B$ can be reconstructed, which are dual to each other.\\

Observe that the composites of dec/rec maps also provides us with a commutative nabla-like diagram in the present situation, in contrast to starting from two dual theories, which does not necessarily imply that the accompanying nabla-like diagram is commutative.\\ 

As a corollary, from the perspective of constructing a ``theory of everything", if one has two theories $A$ and $B$ that satisfy a formal uncertainty principle, hence are presumably dual to each other, it follows there must be some underlying theory containing the ground state $(C_{AB}, N_A, N_B, E_{AB}) = (C_{AB},N_A,N_B)$ we should be able to project down to.\\

\bigskip
\footnotesize
\noindent
\textit{e-mail address}: \texttt{rg.mathematics@gmail.com}.


\begin{thebibliography}{20}
	\bibitem[A]{A} K. Arakawa, \textit{The Grothendieck Construction for $\infty$-categories Fibered over Categorical Patterns}, arXiv:2404.01025[math.CT].
	\bibitem[RG]{RG} R. Gauthier, \textit{Introduction to the Category of Motivic Spectra}, arXiv:2303.09537[math.AG].
	\bibitem[AG]{AG} A. Grothendieck, \textit{Revetements Etales et Groupe Fondamental}, Seminaire de Geometrie Algebrique du Bois Marie 1960/61.
	\bibitem[HB]{HB} S. de Haro, J. Butterfield, \textit{On Symmetry and Duality}, arXiv:1905.05966v2[physics.hist-ph].
	\bibitem[Lu]{Lu} J. Lurie, \textit{Higher Topos Theory}, Annals of Mathematics Studies, 2009, Princeton University Press.
	\bibitem[Lu2]{Lu2} J. Lurie, \textit{Higher Algebra}, people.math.harvard.edu/$\sim$lurie/papers/HA.pdf.
	\bibitem[AMG]{AMG} A. Mazel-Gee, \textit{On the Grothendieck Construction for $\infty$-Categories}, J. of Pure and Applied Algebra, Vol 223, Issue 11, Nov. 2019, pp.4602-4651.
\end{thebibliography}
\end{document}